\documentclass[12pt,reqno]{amsart}

\usepackage{amssymb,latexsym}

\usepackage{enumerate}

\usepackage[french,english]{babel}
\usepackage{amsmath}
\usepackage{graphicx}
\usepackage{amssymb}
\usepackage{bbm}
\usepackage{amsthm,mathtools}
\usepackage{ulem}
\usepackage{geometry}
\usepackage{tikz-cd}
\usepackage{mathrsfs}
\usepackage[colorinlistoftodos]{todonotes}
\usepackage{enumitem}
\usepackage{verbatim}
\usepackage[foot]{amsaddr}
\usepackage{dsfont}
\usepackage{cite}

\makeatletter

\@namedef{subjclassname@2010}{
	
	\textup{2020} Mathematics Subject Classification}

\makeatother
\newtheorem{thm}{Theorem}[section]
\newtheorem*{thm*}{Theorem}

\newtheorem{lem}[thm]{Lemma}

\theoremstyle{definition}

\numberwithin{equation}{section}

\newcommand{\mcf}{\mathcal{F}}

\newcommand{\mcm}{\mathcal{M}}

\usepackage{hyperref}
\hypersetup{hypertex=true,colorlinks=true,linkcolor=blue,anchorcolor=blue,citecolor=blue}
\newcommand{\newabstract}[1]{%
	\par\bigskip
	\csname otherlanguage*\endcsname{#1}%
	\csname captions#1\endcsname
	\item[\hskip\labelsep\scshape\abstractname.]
}

\begin{document}

	\baselineskip=17pt

	\title[Quadratic character sums with multiplicative coefficients]{Quadratic character sums with multiplicative coefficients}

	\author{Zikang Dong\textsuperscript{1}}
    \author{Zhonghua Li\textsuperscript{2}}
	\author{Yutong Song\textsuperscript{2}}
    \author{Shengbo Zhao\textsuperscript{2}}
	\address{1.School of Mathematical Sciences, Soochow University, Suzhou 215006, P. R. China}
	\address{2.School of Mathematical Sciences, Key Laboratory of Intelligent Computing and Applications(Ministry of Education), Tongji University, Shanghai 200092, P. R. China}
   
	\email{zikangdong@gmail.com}
    \email{zhonghua\_li@tongji.edu.cn}
	\email{99yutongsong@gmail.com}
	\email{shengbozhao@hotmail.com}
	\date{\today}
	
	\begin{abstract} 
		In this article, we study extreme values of quadratic character sums with multiplicative coefficients $\sum_{n \le N}f(n)\chi_d(n)$. For a positive number $N$ within a suitable range, we employ the resonance method to establish a conditional $\Omega$-result.
	\end{abstract}
	
\keywords{Quadratic character sums, extreme values, multiplicative functions,  resonance method. }

	\subjclass[2020]{Primary 11L40, 11M06, 11N37.}
	
	\maketitle
	
	\section{Introduction}

Character sums play a crucial role in the study of modern number theory and have rich applications in various aspects. The study of character sums has a long history since 1918, when P\'olya and Vinogradov separately showed the well-known inequality
$$\sum_{n\le N}\chi(n)\ll \sqrt q\log q.$$
Here $q$ is any large prime and  $\chi({\rm mod}\;q)$ is any non-principal Dirichlet character modulo $q$. Under the assumption of the generalized Riemann hypothesis (GRH), Montgomery and Vaughan \cite{MV77} refined this upper bound
$$\sum_{n\le N}\chi(n)\ll \sqrt q\log\log q.$$
Define the maximum of character sums 
$$M(\chi) \coloneqq \max\limits_{t \leq q} \Big| \sum_{n \leq t} \chi(n) \Big|.$$ 
Paley \cite{Pl} showed that
$$
M(\chi) \gg \sqrt{q}\log \log q
$$
holds for infinitely many moduli $q$ and quadratic characters $\chi \, ({\rm{mod}}\ q)$. Later, this result was refined by Bateman and Chowla \cite{BC} in 1950. They demonstrated that there exists a sequence of moduli $q$ and primitive quadratic characters $\chi \, ({\rm{mod}}\ q)$ such that 
$$
M(\chi) \ge \Big( \frac{\mathrm{e}^\gamma}{\pi} + o(1) \Big) \sqrt{q} \log \log q
$$
holds as $q \to \infty,$ where $\gamma$ is the Euler-Mascheroni constant. In 2007, Granville and Soundararajan \cite{GS07} showed under GRH that 
$$
M(\chi) \le \Big( \frac{2\mathrm{e}^\gamma}{\pi} + o(1) \Big) \sqrt{q} \log \log q.
$$
They conjectured the Omega result of Bateman and Chowla should correspond to the true extreme values of $M(\chi)$, namely that
$$
M(\chi) \le \Big( \frac{\mathrm{e}^\gamma}{\pi} + o(1) \Big) \sqrt{q} \log \log q.
$$

Consider the dual maximum of character sums for any $N\le q$:
$$\max_{\substack{\chi \, ({\rm{mod}}\ q)\\\chi \neq \chi_0 }} \Big| \sum_{n \leq N} \chi(n) \Big|.$$
Granville and Soundararajan \cite{GS01} showed several important lower bounds for it as $N$ varies in different rangs according to $q$. See Theorems 4--8 in  \cite{GS01}. The main method used there is  moments of random multiplicative functions (and also of character sums). They also generalized these results to the case of quadratic character sums. See Theorems 9--11 there.
\par
  Another method to study lower bounds for the above maximum is the resonance method. Since Soundararajan \cite{Sound} improved Voronin's resonance method in \cite{Vor}, this approach has been widely applied to the problem of detecting  large values of some important  functions, such as the Riemann zeta function and Dirichlet $L$-functions.  For more details on this topic, we refer to \cite{Ais, AMM, AMMP, BS2, BS, BT, Hou}. In 2012, Hough \cite{Hou} employed this method to investigate large values of character sums. Denote
$$
\Delta(N, q) \coloneqq \max_{\substack{\chi \, ({\rm{mod}}\ q)\\\chi \neq \chi_0 }} \Big| \sum_{n \leq N} \chi(n) \Big|.
$$
Hough established 
$$
\Delta(N, q) \geq \sqrt{N} \exp \bigg( (1 + o(1)) \sqrt{ \frac{(1 - \theta) \log q}{\log \log q} } \bigg),
$$
where $q$ is a large prime and $N=q^{\theta}$ lies within a suitable range. Later, this result was imporved by La Bretèche and Tenenbaum  \cite{BT}. Similar results for smaller $N$ was established by Munsch \cite{Mun}.
\par
Quadratic character sums, as an important special case, are closely related to quadratic fields. Relevant research can be traced back to Gauss's work on binary quadratic forms. Let $\mathcal{D}$ denote the set of all fundamental discriminants, i.e.
$$
\mathcal{D} := \{d\in{\mathbb Z}: d \,\, \text{is a fundamental discriminant} \},
$$
and $\chi_d = \big( \frac{d}{\cdot}\big)$ be the real primitive character associated to the fundamental discriminant $d$.
\par

As we mentioned previously, in 2001, Granville and Soundararajan \cite{GS01} established several lower bounds for the following maximum
$$
\max_{\substack{X<|d|\le2X \\ d\in \mathcal{D}}}  \sum_{m \le M} \chi_d(m) .
$$
 Kalmynin \cite{Kal} gave a strengthening one of their results. Lamzouri \cite{Lam} conducted a relatively comprehensive study on the distribution of large values of quadratic character sums and applied it to the problem of the positivity of sums of the Legendre symbol, which was previously studied by Fekete, Chowla, Montgomery and so on.  See \cite{BC, GS01, Hou, Kal, Lam} for more details. 
\par
  For character sums with multiplicative coefficients $\sum_{n\le N}f(n)\chi(n)$, Harper \cite{Har} showed that if $|f(p)|=1$ and $|f(p^j)|\le1$ for any prime $p$ and integer $j$, then the low moments ($0<k<1$ fixed) can be bounded with better than square-root cancellation:
 $$\frac{1}{\varphi(q)}\sum_{\chi({\rm mod}\;q)}\Big|\sum_{n\le N}f(n)\chi(n)\Big|^{2k}\ll \bigg(\frac{N}{\sqrt{\log_2(\min\{N,q/N\})}}\bigg)^k.$$
Similar result for zeta sums also holds. Harper conjectured that, if $f(n)$ is not identically equal to 1, then the minimum $\min\{N,q/N\}$ in the above inequality can be replaced by $N$. This means that the ``Fourier flips" for character sums no longer influence the mix sums when $\sqrt q<N<q$. Thus the character sums with multiplicative coefficients also play important roles in analytic number theory as character sums. Recently, Xu and Yang \cite{XY} also studied extreme values of zeta sums with multiplicative coefficients
$$
  \Big|\sum_{n \le N}f(n)n^{it} \Big|
$$

 In this paper, motivated by the work of Harper \cite{Har} and Xu and Yang \cite{XY},
 we focus on extreme values of quadratic character sums with multiplicative coefficients
$$
\sum_{n \le N}f(n)\chi_d(n),
$$
where $f(n)$ is a multiplicative function.

Use $\square$ to indicate a perfect square number. Then,  we define the set  $\mcf$ 
\begin{align*}\mcf&:=\\&\big\{f\;{\rm completely\; multiplicative}: \;|f(k)|=1,\;\forall\, k\in{\mathbb N};\;{\rm{Re}}f(n)\overline{f(m})\ge 0,\;\forall \,mn=\square\big\}.\end{align*} Under  GRH, we have the following theorem, which gives conditional extreme values of quadratic character sums with multiplicative coefficients.
\begin{thm}\label{thm1.1}
Fix $\delta \in (0,\frac{1}{100})$. Let $X$ be a sufficiently large number and $N$ be any positive number satisfying $\exp((\log X)^{\frac{1}{2}+\delta}) \le N \le X^{\frac{1}{4}-\delta}$. Assuming GRH is true, then we have
    $$ \max_{\substack{X<|d|\le2X \\ d\in \mathcal{D}}} \Big|\sum_{n\le N}f(n)\chi_d(n)\Big| \ge \sqrt{N}\exp\bigg(\Big(\tfrac{\sqrt{2}}{2}+o(1)\Big){\sqrt{\frac{\log (X/N^4)}{\log\log (X/N^4)}}}\bigg)$$
for all $f \in \mcf$.
\end{thm}
\par
Theorem \ref{thm1.1} shows that, for all functions $f \in \mcf,$ the order of extreme values of such quadratic character sums is similar to that of \cite{XY}. Incidentally, in the author's recent work \cite{DLSZ}, we investigated extreme values of character sums $\sum_{n \le N}f(n)\chi(n)$ for Dirichlet characters $\chi$ modulo a large prime $q$, and showed a similar result  without assuming GRH.
\par
For ease of representation, we now introduce some notations. Let $\varepsilon $ be an arbitrarily small positive number. Note that each appearance of $\varepsilon$ does not necessarily denote the same value. 
We set $\log_j x {\rm{:}}= \log_{j-1} \log x$, the $j$-th iterated logarithm. Finally, let $\mu$ be the Möbius function as usual.

      \section{Preliminary Lemmas}
      In this section, we give some lemmas that we will use later. The first lemma below serves to handle the partial sums of quadratic characters $\chi_d(n)$, which is essentially  \cite[Lemma 1]{DM}.
\begin{lem}
    \label{lemmaDM}
    Assume GRH is true. Let $n=n_1 n_2^2$ be a positive integer with $n_1$ the square-free part of $n$. Then for any $\varepsilon>0$, we have
    $$
    \sum_{\substack{|d| \leq X \\ d \in \mathcal{D}}} \chi_d(n)=\frac{X}{\zeta(2)} \prod_{p \mid n}\Big(\frac{p}{p+1}\Big) \mathrm{1}_{n=\square}+O\Big(X^{\frac{1}{2}+\varepsilon} g_1\left(n_1\right) g_2\left(n_2\right)\Big),
    $$
    where $g_1(n_1)=\exp \big(\left(\log n_1\right)^{1-\varepsilon} \big)$, $g_2\left(n_2\right)=\sum_{q \mid n_2} \mu^2(q) q^{-(1 / 2+\varepsilon)}$, and $\mathrm{1}_{n=\square}$ indicates the indicator function of the square numbers.
\end{lem}
\par
To simplify the writing, we can provide rough upper bounds for $g_1$ and $g_2$. In fact, by the definition of $g_1(n)$, we directly obtain
$$
    g_1(n_1) = \exp \Big( (\log n_1)^{1-\varepsilon}\Big) \ll n_1^\varepsilon \ll n^\varepsilon.
$$
And then, we have
$$
    g_2(n_2) = \sum_{q \mid n_2}\frac{\mu(q)^2}{q^{\frac{1}{2}+\varepsilon}} = \prod_{p \mid n_2}\Big(1+\frac{1}{p^{\frac{1}{2}+\varepsilon}} \Big) \ll n_2^{\varepsilon} \ll n^\varepsilon.
$$

\par
The second lemma states the lower bound for the ratio of two sums concerning a special multiplicative function $r(\cdot)$, which plays a crucial role in the resonance method. This conclusion follows directly from \cite[Eq. (2.8)]{XY}, which is based on the work of \cite{Hou}.
\begin{lem}\label{lem2.2}
    Let $Y$ be large and $\lambda=\sqrt{\log Y\log_2 Y}$. Define the multiplicative function $r$ supported on square-free integers and for any prime $p$:
$$r(p)=\begin{cases}
   \frac{\lambda}{\sqrt p \log p}, &  \lambda\le p\le \exp((\log\lambda)^2),\\
   0, & {\rm otherwise.}
\end{cases}$$
If $\log N>3\lambda\log_2\lambda$, then we have
\begin{equation}\label{DD}
    \sum_{a,b\le Y}\sum_{m,n\le N\atop an=bm}r(a)r(b)\Big/\sum_{n\le Y}r(n)^2\ge N\exp{\bigg((2+o(1))\sqrt{\frac{\log Y}{\log_2 Y}}}\bigg).
\end{equation}
\end{lem}

        \section{Proof of Theorem \ref{thm1.1}}
    We utilize the resonance method proposed by Soundararajan \cite{Sound} to prove Theorem \ref{thm1.1}. For convenience, we begin by introducing some definitions. Let 
    $$
    R_d := \sum_{n \leq Y} r_f(n) \chi_d(n)
    $$
    be the resonator, where $r_f(n)$ is an arithmetical function that will be chosen later and $Y = X^{1/2-\delta}/N^{2}.$ Then, we define the following two sums
    $$
    S_1= S_1(X,R_d) := \sum_{\substack{X<|d| \leq 2 X \\ d \in \mathcal{D}}}  |R_d|^2, 
    $$
    and
    $$
    S_2= S_2(X,R_d) := \sum_{\substack{X<|d| \leq 2 X \\ d \in \mathcal{D}}}  |D_N(f)|^2 |R_d|^2,
    $$
    where
    $$
    D_N(f) = D_{N,d}(f) := \sum_{n \leqslant N} f(n) \chi_d(n).
    $$
    \par
    Plainly, we have 
    \begin{equation}
    \label{maxD}
    \max_{\substack{X<|d| \leq 2X \\ d \in \mathcal{D}}} |D_N(f)| \geq \sqrt{\frac{S_2}{S_1}}.
    \end{equation}
    \par
First, we need an effective lower bound for $S_2$. Plugging the definition of $D_N(f)$ and $R_d$ into $S_2$ yields 
$$
S_2 = \sum_{\substack{X<|d| \leq 2 X \\ d \in \mathcal{D}}} \sum_{m,n\leq N}f(n)\overline{f(m)}\chi_d(mn) \sum_{a,b \leq Y}r_f(a)\overline{r_f(b)}\chi_d(ab).
$$
To avoid the effect of $f \in \mcf$, we set $r_f(n)=f(n)r(n)$, where $r$ is a non-negative multiplicative function to be chosen later. We change the order of summation and obtain
$$
 S_2 = \sum_{a,b\leq Y}\sum_{m,n\le N} f(n)\overline{f(m)}  f(a)r(a)\overline{f(b)}r(b)\sum_{\substack{X<|d| \leq 2 X \\ d \in \mathcal{D}}} \chi_d(abmn).
$$
We split the sum into two parts: the case $abmn = \square$ and the case $abmn \neq \square$. According to Lemma \ref{lemmaDM}, we have
\begin{align*}
    S_2 &= \sum_{a,b\leq Y }\Big( \sum_{\substack{m,n \leq N \\ abmn = \square}} +\sum_{\substack{m,n \leq N \\ abmn \neq \square}} \Big) f(n)\overline{f(m)}  f(a)r(a)\overline{f(b)}r(b) \sum_{\substack{X<|d| \leq 2 X \\ d \in \mathcal{D}}} \chi_d(a b m n) \\
    &= \frac{X}{\zeta(2)}\sum_{a,b\leq Y } \sum_{\substack{m,n \leq N \\ abmn = \square}} f(n)\overline{f(m)} f(a)\overline{f(b)}r(a)r(b)\prod_{p \mid abmn}\Big( \frac{p}{p+1}\Big) \\
    & \quad \,+ O\Big( X^{\frac{1}{2}+\varepsilon}N^{3\varepsilon} Y^{3\varepsilon} \sum_{\substack{ a,b \le Y \\ m,n \le N}}r(a)r(b)\Big).
\end{align*}
The Cauchy-Schwarz inequality yields
\begin{align}\label{S2M}
        S_2 &= \frac{X}{\zeta(2)}\sum_{a,b\leq Y } \sum_{ m,n \leq N\atop abmn=\square}f(n)\overline{f(m)} f(a)\overline{f(b)}r(a)r(b)\prod_{p \mid abmn}\Big( \frac{p}{p+1}\Big) \nonumber \\
         & \quad \, +\, O\Big( X^{\frac{1}{2}+\varepsilon} N^{2+3\varepsilon} Y^{1+3\varepsilon} \sum_{n \le Y}r(n)^2 \Big).
\end{align}
\par
We use $\frac{X}{\zeta(2)} \mcm$ to denote the main term of $S_2$ in \eqref{S2M}. Note that terms in $\mathcal{M}$ satisfying $an=bm$ are independent of the completely multiplicative function $f$. Define the arithmetic function $h(\cdot)$ by $h(1)=1$ and $h(n) = \prod_{p \mid n} \frac{p}{p+1}$ for $n \neq1.$
Then we have
\begin{align}\label{Mlower}
    \mcm & =\sum_{a,b\leq Y }\sum_{\substack{ m,n \le N \\ a n=b m}}r(a)r(b)h(an) + \sum_{a,b\leq Y } \sum_{\substack{ m,n \le N\\abmn=\square \\ a n\neq b m}} f(an)\overline{f(bm)}r(a)r(b)h(abmn) \nonumber \\
    &= \sum_{a,b\leq Y } \sum_{\substack{ m,n \le N \\ a n=b m}}r(a)r(b)h(an) + 2{\rm{Re}}\sum_{a,b\leq Y } \sum_{\substack{ m,n \le N\\abmn=\square \\ a n > b m }} f(an)\overline{f(bm)}r(a)r(b)h(abmn) \nonumber \\
    & \ge   \sum_{a,b\leq Y }\sum_{\substack{ m,n \le N \\ a n=b m}}r(a)r(b)h(an), 
\end{align}
where, in the last step, we employ our assumption $f \in \mcf$ and the fact that $h(m) >0$ and $r(m) \ge 0$ for all $m \in \mathbb{N}$.
Plainly, $h$ is a decreasing function. By the prime number theorem, we have
\begin{equation}
    \label{hlowerbound}
    h(an) \ge \prod_{p \le NY}\frac{p}{p+1} \ge \prod_{p \le X}\frac{p}{p+1}= \exp \bigg(  \sum_{p \le X}\log\Big(1-\frac{1}{p+1}\Big) \bigg) \ge(\log X)^{-c},
\end{equation} for some absolutely positice $c$.
Plugging \eqref{hlowerbound} into \eqref{Mlower}, we obtain
\begin{equation}\label{Mlowerbound}
    \mcm \ge  (\log X)^{-c} \sum_{\substack{a,b\le Y \\ m,n \le N \\ a n=b m}}r(a)r(b).
\end{equation}
Combining \eqref{S2M} and \eqref{Mlowerbound} and employing $Y=X^{1/2-\delta}/N^{2},$ we have
\begin{equation}
    \label{S2}
    S_2 \ge \frac{X}{\zeta(2)} (\log X)^{-c} \sum_{a,b\le Y}\sum_{\substack{ m,n \le N \\ a n=b m}}r(a)r(b) +O\Big(X^{1-\delta+7\varepsilon}\sum_{n \le Y}r(n)^2\Big).
\end{equation}
Since $\varepsilon$ is arbitrarily small, we can request  $\varepsilon < \frac{\delta}{100}$.
%估计S_1的上界
\par
Then we compute the upper bound for $S_1$. Substituting the definition of $R_d$ into $S_1$ and changing the order of summation, we get
$$
S_1 = \sum_{a,b \leq X} f(a)\overline{f(b)}r(a)r(b) \sum_{\substack{X<|d| \leq 2 X \\ d \in \mathcal{D}}} \chi_d(ab). 
$$
Splitting $ S_1 $ into two parts by considering two cases $a=b$ and $a \neq b$ and using the fact that $|f(n)|=1$, we have
$$
S_1 = \sum_{n \leq Y}r(n)^2 \sum_{\substack{X<|d| \leq 2 X \\ d \in \mathcal{D}}} \chi_d(n^2) + \sum_{\substack{a,b\leq Y \\ a \neq b}} f(a)\overline{f(b)}r(a)r(b) \sum_{\substack{X<|d| \leq 2 X \\ d \in \mathcal{D}}} \chi_d(ab).
$$
Then, employing Lemma \ref{lemmaDM} and using the Cauchy-Schwarz inequality yields that
$$
    S_1 = \frac{X}{\zeta(2)}\sum_{n 
    \leq Y} r(n)^2 h(n) + O\Big(X^{\frac{1}{2}+\varepsilon} Y^{1+2\varepsilon} \sum_{n \leq Y}r(n)^2\Big).
$$
Since $h(n) \le 1$, we obtain that
$$
  S_1 \le \Big( \frac{X}{\zeta(2)} + O\big(X^{\frac{1}{2}+\varepsilon} Y^{1+2\varepsilon} \big)\Big)\sum_{n \leq Y} r(n)^2.
$$
So the upper bound for $S_1$ can be established as follows
\begin{equation}
    \label{S1upperbound}
    S_1 \leq \frac{X}{\zeta(2)}(1+o(1))\sum_{n \leq Y} r(n)^2.
\end{equation}
\par
We choose the function $r(\cdot)$ to coincide with that in Lemma \ref{lem2.2}. Subsequently, combining \eqref{S2} and \eqref{S1upperbound}, we have 
\begin{equation*}
        \frac{S_2}{S_1} \gg  (\log X)^{-c}\sum_{a,b\le Y}\sum_{\substack{m,n\le N\\an=bm}}r(a)r(b)\Big/\sum_{n\le Y}r(n)^2 .
\end{equation*}
Employing \eqref{DD} in Lemma \ref{lem2.2} gives 
\begin{equation*}
    \frac{S_2}{S_1} \ge N \exp{\bigg((2+o(1))\sqrt{\frac{\log Y}{\log_2 Y}}}\bigg).
\end{equation*}
Finally, combining \eqref{maxD} and recalling that $Y = X^{1/2-\delta}/N^{2}$, we complete the proof of Theorem \ref{thm1.1}.

	\section*{Acknowledgements}
	The authors are supported by the Shanghai Magnolia Talent Plan Pujiang Project (Grant No. 24PJD140) and the National
	Natural Science Foundation of China (Grant No. 	1240011770).

	\normalem


\begin{thebibliography}{99}
		
			
\bibitem{Ais}
Aistleitner, C. {\emph Lower bounds for the maximum of the Riemann zeta function along vertical lines}, {\it Math. Ann.}, {\bf 365} (2016), 73--96.

\bibitem{AMM}
Aistleitner, C.; Mahatab, K.; Munsch, M. {\emph Extreme values of the Riemann zeta function on the $1$-line}, {\it Int. Math. Res. Not.}, {\bf 22} (2019), 6924--6932.

\bibitem{AMMP}
Aistleitner, C.; Mahatab, K.; Munsch, M.; Peyrot, A. {\emph On large values of $L(\sigma, \chi)$}, {\it Q. J. Math.}, {\bf 70} (2019), 831--848.

\bibitem{BC} Bateman, P. T.; Chowla S. {\emph Averages of character sums}, {\it Proc. Amer. Math. Soc.}, {\bf 1} (1950), 781--787.

\bibitem{BS2} Bondarenko, A.; Seip, K. {\emph Large greatest common divisor sums and extreme values of the
	Riemann zeta function}, {\it Duke Math. J.}, {\bf 166} (2017), 685--701.
	
\bibitem{BS} Bondarenko, A.; Seip, K. {\emph Extreme values of the Riemann zeta function and its argument}, {\it  Math. Ann.}, {\bf 372} (2018), 999--1015.
%\bibitem{Buj} Bujold, C. {\emph Long large character sums}, {\it  Mathematika}, {\bf 68}, (2022), 15--50.

\bibitem{BT} de la Bret\`{e}che, R.; Tenenbaum, G. {\emph Sommes de G\'{a}l et applications. (French) [G\'{a}l-type sums and applications]}, {\it Proc. Lond. Math. Soc.}, {\bf 119} (2019), 104--134.

\bibitem{DM} Darbar, P.; Maiti, G. {\emph Large values of quadratic Dirichlet $L$-functions}, {\it Math. Ann.}, (2025), pp. 1--33.

\bibitem{DLSZ} Dong, Z. ; Li, Z. ; Song, Y. ; Zhao, S. {\emph Large values of character sums with multiplicative coefficients}, {ArXiv:2508.09750.}

\bibitem{GS01} Granville, A.; Soundararajan K. {\emph Large character sums}, {\it J. Amer. Math. Soc.}, {\bf 14} (2001), 365--397.
\bibitem{GS07} Granville, A.; Soundararajan K. {\emph Large character sums: pretentious characters and the P\'olya-Vinogradov theorem}, {\it J. Amer. Math. Soc.}, {\bf 20} (2007), 357--384.

\bibitem{Har} Harper, A. J. \emph{The typical size of character and zeta sums is $o(\sqrt{x})$}, {ArXiv:2301.04390.}

\bibitem{Hou} Hough, B. {\emph The resonance method for large character sums}, {\it  Mathematika},  {\bf 59},(2013), 87--118.

\bibitem{Kal}  Kalmynin, A.B. {\emph Large values of short character sums},  {\it J. Number Theory},  {\bf 198} (2019), 200--210.

\bibitem{Lam} Lamzouri, Y. {\emph The distribution of large quadratic character sums and applications}, {\it  Algebra Number Theory}, {\bf 18} (2024), 2091--2131.

\bibitem{MV77}  Mntgomery H. L.; Vaughan R. C. {\emph Exponential sums with multiplicative coefficients}, {\it Invent. Math}, {\bf 43} (1977), 69--82.

\bibitem{Mun} Munsch, M.  {\emph The maximum size of short character sums}, {\it Ramanujan J.} {\bf 53} (2020), 27–-38.

\bibitem{Pl} Paley, R. E. A. C. {\emph A theorem on characters}, {\it J. London Math. Soc.}, {\bf 7} (1932), 28--32.

\bibitem{Sound}
 Soundararajan, K. {\emph Extreme values of zeta and $L$-functions}, {\it Math. Ann.}, {\bf 342} (2008),
67--86. 

\bibitem{Vor} Voronin, S. M.; {\emph Lower bounds in Riemann zeta-function theory}, {\it Izv. Akad. Nauk SSSR Ser. Mat.}, {\bf 52} (1988), 882--892, 896.
        
\bibitem{XY} Xu, M.W.; Yang, D. {\emph Extreme values of Dirichlet polynomials with multiplicative coefficients}, {\it J. Number Theory}, {\bf 258} (2024), 173--180.

	\end{thebibliography}
\end{document}